\documentstyle[12pt]{amsart} 
\newcommand{\nl}{\newline}

\newcommand{\be}{\begin{enumerate}}
\newcommand{\ee}{\end{enumerate}}
\newcommand{\tbe}{(to be extended) }
\newcommand{\pr}{{\bf Proof. }}

\newcommand{\ca}{{\cal A}}
\newcommand{\cb}{{\cal B}}
\newcommand{\cc}{{\cal C}}
\newcommand{\cd}{{\cal D}}
\newcommand{\ce}{{\cal E}}
\newcommand{\cf}{{\cal F}}
\newcommand{\cm}{{\cal M}}
\newcommand{\cn}{{\cal N}}
\newcommand{\co}{{\cal O}}
\newcommand{\cp}{{\cal P}}
\newcommand{\cs}{{\cal S}}
\newcommand{\ct}{{\cal T}}
\newcommand{\crxc}{{\cal R}_{X,c}}
\newcommand{\bn}{\cb_n}
\newcommand{\tbn}{\hat\cb_n}
\newcommand{\bd}{\textrm{Bd}\,}
\newcommand{\dbd}{{\cd^{\text{Bd}}}}
\newcommand{\dbdn}{\cd_n^{\text{Bd}}}

\newcommand{\bnk}{{\cal B}_{n,k}}
\newcommand{\bnko}{\bnk^o}
\newcommand{\hbnk}{\hat\bnk}
\newcommand{\hbnko}{\hat\bnko}

\newcommand{\bnkt}{{\cal B}_{n,k,t}}
\newcommand{\qd}{$\Box$}
\newcommand{\th}{{\tilde H}}

\newcommand{\dl}{\text{\bf plim}\,}
\newcommand{\holim}{\text{\bf holim}\,}
\newcommand{\da}{\Delta}
\newcommand{\dan}{\da_n}
\newcommand{\zp}{{\Bbb Z}_p}

\newcommand{\la}{\leftarrow}
\newcommand{\lla}{\longleftarrow}
\newcommand{\ra}{\rightarrow}
\newcommand{\dar}{\downarrow}
\newcommand{\lra}{\longrightarrow}
\newcommand{\follows}{\Rightarrow}
\newcommand{\tcdn}{\tilde\cd_n}
\newcommand{\bo}{\partial\,}
\newcommand{\sm}{\setminus}
\newcommand{\ti}{\tilde }
\newcommand{\st}{\text{st}}
\newcommand{\lk}{\text{lk}\,}
\newcommand{\cat}{\text{\bf Cat}}
\newcommand{\tp}{\text{\bf Top}}
\newcommand{\dd}{{\Bbb D}}
\newcommand{\sd}{\text{sd}}
\newcommand{\wsd}{\widetilde{\sd}}
\newcommand{\il}{\lim}
\newcommand{\id}{\text{id}}
\newcommand{\fns}{\footnotesize}
\newcommand{\sr}{\stackrel}
\newcommand{\sh}{\diamond }
\newcommand{\bu}{\bullet}

\newtheorem{thm}{Theorem}[section]
\newtheorem{df}  [thm]{Definition}
\newtheorem{mthm}[thm]{Main Theorem}
\newtheorem{lm}  [thm]{Lemma}
\newtheorem{crl} [thm]{Corollary}
\newtheorem{prop}[thm]{Proposition}

\numberwithin{equation}{section}

  \title{Diagrams of classifying spaces and $k$-fold Boolean algebras}
              
                \author{Eric Babson and Dmitry N. Kozlov}
                
                          \date{\today \\[0.1cm]
 Research at MSRI was supported in part by NSF grant DMS-9022140.  
 The first author was supported by an NSF Postdoctoral Fellowship. 
 The second author was supported 
by the Swedish Science Council Postdoctoral Fellowship M-PD 11292-303}

\address {Mathematical Science Research Institute, Centennial Drive 1000, 
Berkeley, CA, USA}
\email{ babson@@msri.org, kozlov@@msri.org, kozlov@@math.kth.se} 
 
\begin{document}

\begin{abstract}  
   In this paper we study the problem of determining the homology groups of a 
quotient of a topological space by an action of a group. The method is to represent 
the original topological space as a homotopy limit of a diagram, and then act with 
the group on that diagram. Once it is possible to understand what the action of 
the group on every space in the diagram is, and what it does to the morphisms, we 
can compute the homology groups of the homotopy limit of this quotient diagram.   
  
  Our motivating example is the symmetric deleted join of a simplicial complex. It 
can be represented as a diagram of symmetric deleted products. In the case where the 
simplicial complex in question is a simplex, we perform the complete computation of 
the homology groups with $\zp$ coefficients. For the infinite simplex the spaces in 
the quotient diagram 
are classifying spaces of various direct products of symmetric groups and diagram 
morphisms are induced by group homomorphisms. Combining Nakaoka's description of the 
$\zp$-homology of the symmetric group with a spectral sequence, we reduce the 
computation to an essentially combinatorial problem, which we then solve using the 
braid stratification of a sphere. Finally, we give another description of the problem 
in terms of posets and complete the computation for the case of a finite simplex.  
\end{abstract}

\maketitle

                 \section{Introduction}

  Assume that we have a topological space $B$ and a group $G$ acting on it. We wish 
to compute the homology groups of $B/G$. Sometimes it is possible to find a diagram 
$\cd$ such that $B$ is homotopy equivalent to $\holim\cd$ and $G$ acts on $\cd$ in a 
sufficiently simple way that we can understand the spaces and morphisms of the 
quotient diagram $\cd/G$. We can then compute the homology groups of $\holim(\cd/G)$ 
(for example using the Mayer-Vietoris spectral sequence) and obtain the desired answer. 

Our motivating example is the $n$-fold
deleted join of a simplicial complex $C$ and $G=\cs_n$ acting on $B$ by 
permuting the terms.  We then represent $B$ as the homotopy
limit of a diagram $\cd$, where spaces are various deleted products of $C$, 
and subsequently act with $\cs_n$ on $\cd$.

  We perform the complete computation in the case where $C$ is a simplex. The diagram
for the infinite case consists of different Eilenberg-MacLane spaces
(of various direct products of symmetric groups) combined to
obtain the symmetric deleted join of the infinite simplex.

   In order to actually compute the homology groups, we combine the description of 
the homology groups, with coefficients in $\zp$, of the symmetric group given by 
M.~Nakaoka, \cite{N60} with the standard spectral sequence converging to the homotopy 
limit of a diagram (cf. \cite{Da}). This spectral sequence collapses at the second step. 
We compute the tableaux $E_2^{*,*}$ using a geometric argument involving stratification 
of a sphere by hyperplanes from the braid arrangement. 

The topological spaces that we study can also be described as classifying spaces of 
certain posets with a natural combinatorial description. We use this language to 
exhibit the sequence of symmetric deleted joins of simplices and to prove that their homology
strongly converges to that of the computed example. As a byproduct, we obtain a family
of small finite dimensional complexes, whose homology essentially coincides with the
initial segment of the homology of the symmetric group.

  Here is a brief outline of the contents of the paper.

{\bf Section 2.} The general framework of the studied problem is described in 
terms of categories, double categories and colimits.  

{\bf Section 3.} We describe how these abstract results specialize to our context.

{\bf Section 4.} We explain how to decompose a (symmetric) deleted join of a 
sim\-pli\-cial complex $C$ as a diagram consisting of (symmetric) deleted products 
of~$C$.

{\bf Section 5.} In this long section we compute the homology groups of the
symmetric deleted join of a simplex.

{\bf Section 6.} We interpret our results in terms of posets ($k$-fold Boolean algebras)
and the finite dimensional approximations.

                 \section{The general framework}

 In this section we describe the general framework of our studies in terms of category
theory. We refer to \cite{Ma} as the general reference and [B\'e,Er,Pa] for the 
information concerning double categories.  
Some of the lemmata are stated without proof and should be straightforward to check.  

 When $K$ is a category, we denote by $\co(K)$ the set of its objects and by $\cm(K)$ 
the set of its morphisms.  Sometimes $\co$ will be considered as the set 
of identities in $\cm$.  
If $a,b\in\co(K)$, we denote by $\cm_K(a,b)$ the set of all
morphisms from $a$ to $b$. If $m\in\cm_K(a,b)$, we write $\bo^\bu m=a$ and $\bo_\bu m=b$. 
When $K_1$ and $K_2$ are two categories we denote by $\cf(K_1,K_2)$ the set of all 
functors between $K_1$ and $K_2$. Finally, we let $\ct(F_1,F_2)$ denote the set of all 
natural transformations between two functors $F_1,F_2\in\cf(K_1,K_2)$.

In the following, $\cp$ denotes the category of all posets, $\tp$ - the category of all 
topological spaces, $\cat$ - the category of all categories and {\bf 1} the terminal 
category with one object and only the identity morphism. $\da$ denotes the realization 
functor $\da:\cat\ra\tp$ mapping a category to its realization space, see \cite{Q73,Q78,Se}.

\begin{df} \label{didf}
 Let $S$ and $R$ be categories such that objects of $R$ have a name {\bf Name}.
 A functor $F:S\ra R$ is called a {\bf diagram} of {\bf Name} over $S$.
\end{df}

 For example, when $S$ is a poset and $R=\tp$ we talk about a diagram of to\-po\-logical
spaces over a poset. If $\cd=(F:S\lra R)$ is a diagram we write 
$F=F^\cd$, $S=S^\cd$ and $R=R^\cd$. When it does not lead to confusion,
we simply use $\cd$ instead of $F^\cd$. 

\begin{df}\label{dddf}
  Define the functor $\dd:\cat\ra\cat$. If $K\in\co(\cat)$, $\dd K$ is the category 
of all diagrams of objects of $K$. More specifically, $\dd K$ is defined by:
\begin{itemize}
\item $\co(\dd K)=\{\cd\,| R^\cd=K\}$;
\item $\cm_{\dd K}(A,B)=\{(f,\tau)\,|\,f\in\cf(S^A,S^B) \text{ and }\,\tau\in\ct(F^A,F^B\circ f)\}$. 
\end{itemize}
Furthermore, for any functor $F:K\ra L$, $\dd F$ is a functor from $\dd K$ to
$\dd L$ defined by $\dd F(A)=F\circ A$ and $\dd F((f,\tau))=(f,F \circ \tau)$.
\end{df}

 When $T=(f,\tau)\in\cm_{\dd K}(A,B)$, write $f=f^T$, $\tau=\tau^T$.

Define the functor $\Sigma\in\cf(\dd K,\cat)$ by $\Sigma(X)=S^X$. 

\begin{df}
  Let $X\in\co(\dd\cat)$, $X=(F:S\ra \cat)$. A {\bf sink} of $X$ is $L\in\co(\cat)$ together 
with a collection of morphisms $\{\lambda^L_s\in\cm_{\cat}(F(s),L)\}_{s\in\co(S)}$, such that 
if $\alpha\in\cm_S(s_1,s_2)$ then $\lambda^L_{s_2}\circ\alpha=\lambda^L_{s_1}$. When $L$ 
is universal with respect to this property we call it the {\bf colimit} of $X$ and write 
$L=\il X$.
\end{df}

\begin{df} \label{actp}
  Let $K$ be a category and $G$ a finite group. We say that $G$ {\bf acts on} $K$ if 
there is a diagram $X=(F:G\ra\cat)$, such that the object of $G$ maps to $K$. For brevity 
we shall identify $F(g)$ with $g$, for $g\in G$. The quotient category $K/G$ is the 
colimit of $X$ (it always exists since $\cat$ is cocomplete). 
\end{df}

 Let $X\in\co(\dd\cat)$. We define an equivalence relation $\sim$ on the set of morphisms 
$\coprod_{\sigma\in\co(S^X)}\cm(F^X(\sigma))$ by saying that $n\sim n'$ if there exists 
$m\in\cm(S^X)$, such that $F^X(m)(n)=n'$.  Square brackets $[$ $]$ or $[$ $]_X$ will denote the
$\sim$-equivalence classes. 

The following condition will be crucial for our applications.

{\bf Condition A.} If $X\in\co(\dd\cat)$, $\bo^\bu b=\bo^\bu b'$ and $b'\in[b]$ then $b'=b$.

\begin{lm} \label{l3}
Let $X\in\co(\dd\cat)$ and assume that $X$ satisfies condition A. Then $\il X$ can be 
explicitly described in terms of equivalence classes:
$\cm(\il X)=\{[f]\,|\,f\in\coprod_{\sigma\in\co(S^X)}\cm(F^X(\sigma))\}$. Clearly $[f]\in\cm([\bo^\bu f],[\bo_\bu f])$.
\end{lm}

When $A$ is a category and $a\in\co(A)$ we denote by $(a\dar A)$ the {\bf category
of all objects under $a$}, see \cite{Ma}. $(a\dar A)$ is defined by: 
\begin{itemize} 
\item $\co(a\dar A)=\{f\,|\,f\in\cm(A)\text{ and } \bo^\bu f=a\}$;
\item $\cm_{(a\dar A)}(f_1,f_2)=\{h\,|\,h\in\cm_A(\bo_\bu f_1,\bo_\bu f_2) \text{ and }\,h\circ f_1=f_2\}$. 
\end{itemize}

  We define a functor $\wsd\in\cf(\cat,\cat)$ by: $\wsd(S)=S'$, where $\co(S')=\cm(S)$ 
and $\cm_{S'}(f,g)=\{(h,k)\,|\,h,k\in\cm(S),\,k\circ g\circ h=f\}$. We call $\wsd(S)$ 
the {\bf morphism category} of $S$. 
 
\begin{df}
  Let us define the {\bf subdivision functor} $\sd:\dd K\ra\dd(\cat\times K)$.

(1) Take $A\in\co(\dd K)$. Set $S^{\sd A}=\wsd(S^A)$.  For objects set 
$F^{\sd A}(f)=(\bo_\bu f \dar S^A) \times F^A(\bo^\bu f)$, and for morphisms 
$F^{\sd A}((h,k))=(h\circ)\times F(k)$, where 
$h\circ$ is the contravariant functor acting by composition with $h$. 

(2) If $T \in \cm_{\dd K}(A,B)$ then $f^{\sd T}((h,k))=(f^T(h),f^T(k))$ 
and $\tau^{\sd T}_m=\ti f \times \tau^T_{\bo^\bu m}$ where 
%
%
$\ti f$ is a functor from $(b\dar S)$ to $(f(b)\dar R)$ induced by $f$.
\end{df}

\begin{lm}\label{l2}
  Assume that $X\in\co(\dd\cat)$ and $X$ satisfies condition A, then $\dd\wsd(X)$
satisfies condition A.
\end{lm}

\begin{df}
Let $K$ be a category with products and colimits and $\cn\in \cf(\cat,K)$ a functor.  
$\cn\lim$ is the functor from $\dd K$ to $K$ defined by the following composition:
$$\dd K\stackrel{\sd}{\lra}\dd(\cat\times K)\stackrel{\dd(\cn\times\id)}{\lra}\dd K
\stackrel{\il}{\lra}K.$$
\end{df}

Note that if $\cn$ is the trivial functor to the terminal object in $K$ then 
$\cn\lim=\il$.  

\begin{df} \label{ldf}
If $\cn =\id$, resp. $\cn =\da$ (then $K=\tp$), we call $\cn\lim$ the {\bf poset limit}, resp. 
{\bf homotopy limit}, and use the notation $\cn\lim=\dl$, resp. $\holim$.
\end{df}

 When $\cd$ is a diagram of posets over a poset, the first part of the Definition \ref{ldf} specializes 
to the one given in \cite{WZZ}. Namely, the set of the elements of $\dl\cd$ is 
$\{(p,e)\,|\,p\in S^\cd,$ $e\in F^\cd(p)\}$ and the partial ordering is defined by:
$(p,e)\geq (p',e')$ iff $p\geq p'$ and $f_{pp'}(e)\geq$~$e'$ where $f_{pp'}=F^\cd(p \ra p')$.

\begin{prop}
  Assume that we are given two categories $K$ and $K'$ and functors $\cn\in\cf(\cat,K)$, 
$\cn'\in\cf(\cat,K')$ and $W\in\cf(K,K')$, such that $W\circ \cn =\cn '$ and $W$ preserves 
colimits and products. Then the following diagram commutes
$$\begin{array}{ccccccc}
\dd K&\sr{\sd}{\lra}&\dd(\cat\times K)&\sr{\dd(\cn\times\id)}{\lra}&\dd K&\sr{\il}{\lra}&K\\[.1in]
\textrm{\kern-0.3cm\fns$\dd W$}\dar\kern0.3cm&\quad\textrm{\fns(I)}\quad&\kern0.5cm\dar\textrm
{\fns$\dd(\id\times W)$}\kern-0.5cm&\quad\textrm{\fns(II)}\quad&\kern0.3cm\dar\textrm{\fns$\dd 
W$}\kern-0.3cm&\quad\textrm{\fns(III)}\quad&\kern0.17cm\dar\textrm{\fns$W$}\kern-0.17cm\\[.1in]
\dd K'&\sr{\sd}{\lra}&\dd(\cat\times K')&\sr{\dd(\cn '\times\id)}{\lra}&\dd K'&\sr{\il}{\lra}&K'
\end{array}$$
\end{prop}

\pr The commutativity of the first square follows from the definition of $\sd$.
Square (II) commutes, because $W$ preserves products. Square (III) commutes, since
$W$ preserves colimits. 
\qed

\begin{crl} \label{2.4}
 For every diagram $\cd=(F:S\ra \cat)$ we have $$\holim(\dd \Delta(\cd))\simeq\Delta(\dl\cd).$$ 
\end{crl}

\pr Since $\da$ has a weak homotopy inverse, see \cite{FL}, it preserves colimits up to homotopy.  
\qed

{\bf Note.} If $\cd$ satisfies condition A, the last square commutes exactly, so we 
get equality: $\holim(\dd \Delta(\cd))=\Delta(\dl\cd)$.

 {\bf Note.} In the special case when $S$ is a poset and $K=\cp$ see Simplicial Model 
Lemma, \cite[Proposition 3.19]{WZZ}.

 For an arbitrary category $K$, the {\bf barycentric subdivision} of $K$ is the face 
poset of $\da(K)$. We denote it Bd$\,K$. Bd$\,K$ consists of chains of nonidentity morphisms in $K$ ordered
by composition. This notion can be generalized to diagrams.

\begin{df} \label{bddf}
  Given a diagram $\cd \in \co \dd K$, the {\bf barycentric subdivision} of $\cd$ is a 
diagram, which we will denote $\dbd$. Set $S^{\dbd}=\bd S^\cd$.  For 
$\alpha=(x_1\sr{m_1}{\lla}\dots$ $\sr{m_{t-1}}{\lla}x_t)\in\co(\bd S^\cd)$ with  
$x_i\in\co(S^\cd)$ and  
$m_i\in\cm_{S^\cd}(x_{i+1},x_i)$, we define $ F^{\dbd(\alpha)}=F^\cd(x_t)$. For all pairs 
$\alpha,\beta\in\co(\bd S)$, where $\alpha$ is as above and $\beta<\alpha$, 
i.e. $\beta=(x_{i_1}\la\dots\la x_{i_v})$ for some $1\leq i_1<\dots<i_v\leq t$,
the set $\cm_{S^\dbd}(\alpha,\beta)$ has exactly one element $t$ and 
$F^\dbd(t)$ is $\id$ 
if $x_t=x_{i_v}$ and $m_{i_v}\circ m_{i_{v+1}}\circ\dots\circ m_{t-1}$ otherwise.
\end{df}

\begin{df}
  Let $K$ be some given category. $\dd^2 K$ is a category defined by: 
\begin{itemize}
\item $\co(\dd^2 K)=\{(C,f)\,|\,C\text{ is a double category, }f\in\cf(C,K^2)\}$; 
\item $\cm_{\dd^2 K}((C_1,f_1),(C_2,f_2))=\{(f,\tau)\,|\,f\in\cf(C_1,C_2),\,
\tau\in\ct(f_1,f_2\circ f)\}$.
\end{itemize}
Here $K^2$ is the double category whose morphisms are the commuting squares 
in $K$.

 Although this provides a general definition of $\dd^2 K$, for our purposes we need
to assume that conditions 1$\,$-$\,$3 below are satisfied. The morphisms of $C$ can
be composed in two different ways, which we refer to as $\sh$ and $*$ operations.
\end{df}

 Before we formulate the conditions we will need some additional notations.
$C^*$, resp. $C^\sh$, denotes the category which is obtained from the double
category $C$ by taking the same set of morphisms, but considering only $\sh$, 
resp. $*$, operations. $[C]^*$, resp. $[C]^\sh$, is the set of all $*$-, resp. 
$\sh$-connected components. If $a\in\cm(C)$, $[a]^*$, resp. $[a]^\sh$, is the
$*$-, resp. $\sh$-connected component containing $a$. For $a\in\cm(C)$, $\bo^*(a)$,
$\bo_*(a)$, $\bo^\sh(a)$, $\bo_\sh(a)$ will denote upper and lower $*$-identities, 
resp. left and right $\sh$-identities. $C_1^*$ denotes the subcategory of $C^*$
consisting of $\sh$-identities. If $f$ is a double category functor $f:C\ra K^2$, $f^*$, 
$f_*$, $f^\sh$ and $f_\sh$ denote the restrictions of $f$ to the left, right, upper and 
lower morphisms. Clearly, $f^*,f_*\in\cf(C^*,K)$, $f^\sh,f_\sh\in\cf(C^\sh,K)$. 
 
{\bf Condition 1.} If $a,b\in\cm(C)$, $a$ is a $\sh$-identity and $a*b$ exist, then
$a*b$ is a $\sh$-identity again. The same is true if we swap $*$ and $\sh$. 

{\bf Condition 2.} If $a,b,c\in \cm C$ and $a*c$, $a\sh b$ exist, then there 
exists $d\in \cm C$, such that $c\sh d$ and $b*d$ exist. The same is true if we 
swap $*$ and $\sh$ or the order of the first operand. The uniqueness of $d$ follows from 
 Condition 3.


{\bf Condition 3.} If $b'\in[b]^\sh \subseteq \cm C$ and $\bo^\bu b=\bo^\bu b'$  
then $b=b'$ and the same is true if we swap $*$ and $\sh$.

  We shall now define a functor $\varphi_\sh:\dd^2 K\ra\dd\dd K$.

(1) Let $(C,f)\in\co(\dd^2 K)$. Define a new category $S$ by taking $\cm(S)=[C]^\sh$. 
For $A,B\in\cm(S)$, we call $C$ the composition of $A$ and $B$ if there exist $a\in A$, 
$b\in B$ and $c\in C$, such that $a*b=c$. In this case, we simply write $C=A*B$. 

Let us check that $S$ is well defined as a category. Assume that $a'\in A$, $b'\in B$ 
and $c'=a'*b'$. We need to prove that $c'\in C$. Take $a_1\in A$, such that $a\sh a_1$ 
and $a_1\sh a'$ exist. According to the condition 2 there exist $b_1,b_1'\in B$, such 
that $b\sh b_1$, $b_1\sh b_1'$, $a*b_1$, and $a'*b_1'$ exist. Since $b_1'$, $b' \in [b]^\sh$
 we have $b_1'=b'$ by condition 3. Take $c_1=a_1*b_1$,
then $c\sh c_1$ and $c_1\sh c'$ exist and hence $c'\in C$.

  In particular, objects of $S$ are just $\sh$-connected components formed by 
$*$-identities (according to condition 1, a $\sh$-connected component which has
a single $*$-identity, must entirely consist of $*$-identities).

  Observe that, according to conditions 2 and 3, whenever $[a]^\sh *[b]^\sh$ exists, 
there exists a unique $c\in[b]^\sh$, such that $a*c$ exists. We denote this $c$ by 
$\gamma_{[b]^\sh}(d)$. Clearly $\gamma_{[c^\sh]}:[\bo^* c]^\sh\ra[c]^\sh$.  

  Let us now define $F\in\cf(S,\dd K)$. Each $A\in\co(S)$ is a 
$\sh$-subcategory of $C^\sh$. Define $F(A)=f^\sh|_A$.

   Furthermore, take $B\in\cm_S(A,A')$, (which means that $A*B$ and $B*A'$ exist). Let us 
define $F(B)=(\ti B,\tau)$, where $\ti B\in\cf(A,A')$ (considering $A$ and $A'$ 
as subcategories of $C^\sh$) and $\tau\in\ct(F(A),F(A')\circ
\ti B)$. Let $a\in A$, we define $\ti B(a)=\gamma_{A'}(\gamma_B(a))$. By condition 2, 
$\ti B$ is a functor from $A$ to $A'$. For $a\in A$ and $b\in B$, such that $a*b$ 
exists, we define $\tau_a=f^*(b)$. Clearly $\tau\in\ct(F(A),F(A')\circ\ti B)$.

(2) Let $(f,\tau)\in\cm_{\dd^2 K}((C_1,f_1),(C_2,f_2))$. Define $\varphi_\sh(f,\tau)$ to be 
the morphism between $\varphi_\sh(C_1,f_1)$ and $\varphi_\sh(C_2,f_2)$ as follows.
Let $\varphi_\sh(C_1,f_1)=(F_1:S_1\ra\dd K)$, $\varphi_\sh(C_2,f_2)=(F_2:S_2\ra\dd K)$.
$f$ induces a map from $S_1$ to $S_2$ defined by $\ti f:[a]^\sh\ra[fa]^\sh$. This map 
is well defined and is a functor, since $f\in\cf(C_1,C_2)$. Furthermore, define
$\ti\tau\in\ct(F_1,F_2\circ\ti f)$ by letting $\ti\tau_{[a]^\sh}$ be a restriction of 
$(\bar f,\bar\tau)$ to the connected component $[a]^\sh$, where $\bar f$ and $\bar\tau$ 
are $f$ and $\tau$ considered as functor, resp. natural transformation of corresponding
$*$-categories.

  The functor $\phi_*$ is defined analogously.

\begin{lm} \label{l1}
  If $X\in\co(\dd^2 K)$, then $\dd\Sigma(\phi_\sh X)$ satisfies condition A.
\end{lm}

\begin{lm} \label{l4}
  If $L\in\cf(K,K')$ and $X\in\co(\dd K)$, then $\il \dd L(X)=L(\il X)$ iff there is 
some sink $k$ for $X$ satisfying $L(k)=\il\dd L(X)$. 
\end{lm}

\begin{df}
  Take $X\in\co(\dd\dd K)$, $Y$ a sink for $X$ and $y\in\co(S^Y)$. $X_y\in\co(\dd K)$
is defined by: 
\begin{itemize}
\item $\co(S^{X_y})=\coprod_{\sigma\in\co(S^X)}\{u\in\co(F
^X(\sigma))\,|\,
\lambda^Y_\sigma(u)=y\}$; 
\item $\cm_{S^{X_y}}(u,u')=\{m\in\cm_{S^X}\,|\,F^X(m)(u)=u'\}$.
\item $F^{X_y}(u)=F^{F^X(\sigma)}(u)$.
\end{itemize}
\end{df}

\begin{lm} \label{l5}
  If $X\in\co(\dd\dd K)$, $\dd\Sigma(X)$ satisfies condition A, and $Y\in\co(\dd K)$ is 
a sink for $X$, then $Y=\il X$ iff
\begin{enumerate}
\item[($\alpha$)] $\Sigma(Y)=\il(\dd\Sigma(X))$;
\item[($\beta$)] $\il X_y=F^Y(y)$, for all $y\in\co(S^Y)$.
\end{enumerate}
\end{lm}

\begin{lm} \label{l6}
  If $X\in\co(\dd\dd K)$, $\dd\Sigma(X)$ satisfies condition A, $\sigma\in\co(S^X)$, 
and $m\in\cm(F^X(\sigma))$, then $(\dd\sd X)_{[m]}=(\bo_\bu m\dar S^{F^X(\sigma)})\times 
X_{[\bo^\bu m]}$. 
\end{lm}

\begin{thm}\label{cthm}
Let $K$ be a cocomplete category and $\cn\in\cf(\cat,K)$. Assume that if $X\in\co(\dd K)$ and 
$a\in\co(K)$, then $a\times\il X=\il(a\times X)$, where $a\times X\in\co(\dd K)$ is defined
by $F^{a\times X}(m)=\id_a\times F^X(m)$.

Then squares (I), (III) and (I+II) of the following diagram commute. 
$$\begin{array}{ccccccccccc}
\dd^2 K&\kern-0.3cm\sr{\varphi_\sh}{\lra}&\kern-0.3cm\dd\dd K&\kern-0.3cm\sr{\dd\,\sd}{\lra}&
\kern-0.3cm\dd\dd(\cat\times K)&\kern-0.3cm\sr{\dd\,\dd(\cn\times\id)}{\lra}&\kern-0.3cm
\dd\dd(K\times K)&\kern-0.3cm\sr{\dd\mu}{\lra}&\kern-0.3cm\dd\dd K&\kern-0.2cm
\sr{\dd\il}{\lra}&\kern-0.2cm\dd K\\[.1in]

\textrm{\kern-0.1cm\fns$\varphi_*$}\kern-0.1cm\dar\kern0.2cm&\kern-0.3cm\textrm{\fns(I)}&
\kern-0.3cm\kern0.2cm\dar\textrm{\fns$\il$}\kern-0.2cm&\kern-0.3cm&\kern-0.3cm&\kern-0.3cm
\textrm{\fns(II)}&\kern-0.3cm&\kern-0.3cm&\kern-0.3cm
\kern-0.28cm\textrm{\fns$\il$}\dar\kern0.28cm&\kern-0.2cm\textrm{\fns(III)}&\kern-0.2cm
\kern0.2cm\dar\textrm{\fns$\il$}\kern-0.2cm\\[.1in]

\dd\dd K&\kern-0.3cm\sr{\dd\il}{\lra}&\kern-0.3cm\dd K&\kern-0.3cm\sr{\sd}{\lra}&\kern-0.3cm
\dd(\cat\times K)&\kern-0.3cm\sr{\dd\,(\cn\times\id)}{\lra}&\kern-0.3cm\dd(K\times K)&
\kern-0.3cm\sr{\mu}{\lra}&\kern-0.3cm\dd K&\kern-0.2cm\sr{\il}{\lra}&\kern-0.2cm K
\end{array}$$
\end{thm}

\pr Let us first check that the square (I) commutes. Let $(C,f)\in\co(\dd^2 K)$.
This follows from Lemma \ref{l5} since $\dd \il \phi_*(C,f)$ is a sink for 
$\phi_\sh(C,f)$.  Condition $\alpha$ follows from Lemma \ref{l3} giving 
$S^{\phi_*(C,f)}=S^{\il \phi_\sh (C,f)}=[C]^*$.  Condition $\beta$ follows from the 
observation that $(\phi_\sh(C,f))_{[a]^*}=F^{\phi_*(C,f)}([a]^*)$.

%
%
%
%
%
%
%
%
%
  Let us now show that square (I+II) commutes. Let $H=\mu\circ(F\times\id)\circ\sd$. Pick 
$X\in$Im$\,\phi_\sh$, $X\in\co(\dd\dd K)$. According to Lemma \ref{l1}, $\dd\Sigma(X)$
satisfies condition A. According to Lemma \ref{l2}, $\dd\Sigma(\dd H(X))$ satisfies
condition A. According to Lemma \ref{l3}, $\il X$ and $\il\dd H(X)$ have explicit
descriptions in terms of the equivalence classes.
 
   We want to show that $\il\dd H(X)=H(\il X)$. For that we need to check conditions
$(\alpha)$ and $(\beta)$ of Lemma \ref{l5}.

{\bf Check of $(\alpha)$.} Follows from the explicit description of the colimit.
Namely, the objects of $\Sigma(H(\il X))$ are equivalence classes $[m]$ of morphisms
$m\in\cm(\Sigma(F^X(\sigma)))$, for $\sigma\in\co(S^X)$.

{\bf Check of $(\beta)$.} $(H\il X)([m])=\cn(\bo_\bu m\dar S^{F^X(\sigma)})\times 
F^{\il X}([\bo^\bu m])=\cn(\bo_\bu m\dar S^{F^X(\sigma)})\times\il X_{[\bo^\bu m]}$, where 
the first equality follows from the definition of $H$ and the second from Lemma \ref{l5}.

$(\il\dd H(X))([m])=\il(\dd H(X)_{[m]})=\il(\cn(\bo_\bu m\dar S^{F^X(\sigma)}\times 
X_{[\bo^\bu m]})=\cn(\bo_\bu m\dar S^{F^X(\sigma)})\times\il X_{[\bo^\bu m]}$, where
the first equality follows from Lemma \ref{l5}, second from Lemma \ref{l6} and third
from Lemma \ref{l2}.  

Square (III) commutes in general, see \cite[IX 8]{Ma}, but can also be checked using Lemma 
\ref{l5}.  \qed

            \section{Diagrams of posets and the group action}

  A partially ordered set $P$ can be viewed as a category. The objects of 
that category are just the elements of the poset and the morphism $f:\{x\}\lra\{y\}$ 
between the objects $x,y\in P$ exists and is unique iff $x>y$.  Categories 
arising in this way are characterized by the property $|\cm_P(x,y)|+|\cm_P(y,x)|\leq 1$ 
for all $x,y\in\co(P)$. 

  A {\bf poset map} is a functor between two posets viewed as categories, or equivalently an order preserving map.  

  The diagrams which will turn out to be important in this paper, are {\bf diagrams 
of posets}. Rephrasing Definition \ref{didf} a diagram of posets $\cd$ 
over a posset $P$ is a collection of posets $\{\cd(x)\,|\,x\in P\}$ and poset maps 
$\{f_{xy}:\cd(x)\lra\cd(y)\,|\,x,y\in P, x\geq y\}$ 
such that $f_{yz}\circ f_{xy}=f_{xz}$.  In the previous notation 
$S^\cd=P$, $R^\cd=\cp$, $F^\cd(x)=\cd(x)$ and $F^\cd(x\ra y)=f_{xy}$ 
(i.e. in this case the category $K$ is a category 
of all posets with poset maps serving as morphisms). 

{\bf Note.} In this situation the elements of $\dl \cd$ are pairs 
$(x,e)$ with $x \in P$ and $e \in \cd(x)$.  

\begin{prop}
  $\da(\dl\dbd)\simeq\da(\dl\cd)$.
\end{prop}

\pr Let $A=\dl\dbd$, $B=\dl\cd$ and $P=S^\cd$. We shall show that posets $A$ and $B$ are 
homotopy equivalent. Let us describe a map $\phi:A\lra B$. Pick 
$\alpha=x_1<\dots<x_t \in P^{\text{Bd}}$ and $e\in\dbd(\alpha)$. By Definition \ref{bddf} $\,\dbd(\alpha)$ is 
a copy of $\cd(x_t)$ so we can define $\phi((\alpha,e))=(x_t,e)$.

  To verify that $\phi$ is actually a poset map take $(\alpha,e),(\alpha',e')\in A$, 
such that $(\alpha,e)\geq (\alpha',e')$. Assume $\alpha=x_1<\dots<x_t$ and $e\in\cd(x_t)=
\dbd(\alpha)$. Since $\alpha'\leq\alpha$ we can write $\alpha'=x_{i_1}<\dots<x_{i_m}$, 
where $i_1<\dots<i_m\leq t$, clearly $f_{x_t x_{i_m}}(e)\geq e'$. We have 
$\phi(\alpha,e)=(x_t,e)$ and $\phi(\alpha',e')=(x_{i_m},e')$, hence from what is said 
above it follows that $\phi(\alpha,e)\geq\phi(\alpha',e')$.

Let us now show that $\phi$ satisfies Quillen's conditions, \cite{Q78}, and hence 
induces a homotopy equivalence. Let $x\in P$, $e\in\cd(x)$ and consider the poset 
$\phi^{-1}(B_{\leq(x,e)})$. Let $\alpha \in P^{\text{Bd}}$ be a chain in $P$ consisting of a single 
element $x$. Clearly $e\in\dbd(\alpha)$ and 
$\phi((\alpha,e))=(x,e)$. We claim that there always exists a join of $(\alpha,e)$ 
with any other element of $\phi^{-1}(B_{\leq(x,e)})$. In this case 
$\phi^{-1}(B_{\leq(x,e)})$ is join-contractible, which implies that it is contractible.

 If $(\alpha',e')\in\phi^{-1}(B_{\leq(x,e)})$ with $\alpha'=x_1<\dots<x_t$ 
then $x_t\leq x$. Further $(\alpha,e)\vee(\alpha',e')=(\beta,e)$, where $\beta=(x_1<\dots<x_t<x)$
if $x_t\neq x$ and $\beta=\alpha'$ otherwise. Indeed, $(\beta,e)\geq(\alpha,e)$ and 
$(\beta,e)\geq(\alpha',e')$, since $\beta\geq\alpha$, $\beta\geq\alpha'$ and $f_{x x_t}(e)
\geq e'$. Moreover assume that for some $(\tilde\beta,\tilde e)\in\phi^{-1}(B_{\leq(x,e)})$ 
we have that $(\tilde\beta,\tilde e)\geq (\alpha,e)$ and $(\tilde\beta,\tilde e)\geq 
(\alpha',e')$. Then $\tilde\beta$ must include elements $x_1,\dots,x_t$ as well as $x$,
hence $\tilde\beta\geq\beta$. Let $y$ be the maximal element from the chain $\beta$, then 
$y\geq x$ and $f_{yx}(\tilde e)\geq e$ (since $(\tilde\beta,\tilde e)\geq (\alpha,e)$), 
so we can conclude that $(\tilde\beta,\tilde e)\geq(\beta,e)$.
\qed

  In what follows we assume that $P/G$ is a poset and that the conditions of 
Lemma \ref{l3} hold.

  There is a natural poset map $p:P\lra P/G$ which maps $x\in P$ to an orbit 
$a\in P/G$ such that $x\in a$. Obviously every chain in $P$ is mapped to some
chain in $P/G$. As the following lemma shows, the converse is true as well. 

\begin{lm} \label{pre}
  For every chain $a=(a_1<\dots<a_t)$ in $P/G$ there exists at least one chain
$x=(x_1<\dots<x_t)$ in $P$ such that $p(x)=a$.
\end{lm}

\pr We induct on $t$. If $t=1$ then there is nothing to prove. 
If $t=2$ then the statement is a direct consequence of Lemma \ref{l3}.

  Assume $t\geq 3$ and choose a chain $a=(a_1,\dots,a_t)$ with $a_i\in P/G$. 
By induction there is a chain $x_1<\dots<x_{t-1}$ with $x_i\in P$ such that $p(x_i)=a_i$.
According to Lemma \ref{l3}, we can choose $y<z \in P$ such that 
$p(y)=a_{t-1}$, $p(z)=a_t$. Since $p(y)=p(x_{t-1})$ there must exist $g\in G$
such that $g(y)=x_{t-1}$. Put $x_t=g(z)$. Since the action of $G$ is 
order-preserving, we have $y<z\follows g(y)<g(z)\follows x_{t-1}<x_t$. 
Hence, we have found a chain $x_1<\dots<x_t$ with the required properties.
\qed

   The action of $G$ on the poset $P$ naturally induces an action of $G$
on the order simplicial complex $\Delta(P)$ as a topological space, namely
the action on the vertices is given and we extend it over every simplex 
linearly. We denote the quotient space by $\Delta(P)/G$. The map $p$ described 
above induces a map $p^*:\Delta(P)\lra\Delta(P)/G$.

In general $\Delta(P)/G$ will not be a simplicial complex but some simplicial poset, 
see \cite{St}, with every face being a simplex spanned by its vertices, but  
some sets of vertices spanning more than one simplex. 
The following lemma shows that under some 
conditions we obtain a simplicial complex.

\begin{prop}
If  $G$ acts on $P$ in the sense
of the Definition \ref{actp} then the action of $G$ on $P$ induces an action on the
barycentric subdivision $\bd P$ and this new action satisfies the separability 
condition A.
\end{prop}

\pr Let $Q=\bd P$. Choose two chains $x=(x_1,\dots,x_t)$ and $y=(y_1,\dots,y_t)$ in 
$Q$ such that $p(x)=p(y)$. Thus $x_t=(\chi_1,\dots,\chi_m)$ and 
$y_t=(\gamma_1,\dots,\gamma_m)$ are chains in $P$ with $p(x_t)=p(y_t)$. 
Choose $g\in G$ such that $g(x_t)=y_t$.  Since $x_i$ is 
a subchain of $(\chi_1,\dots,\chi_m)$, $y_i$ is a subchain of
$(\gamma_1,\dots,\gamma_m)$, and the action of $g$ 
preserves order, we have $g(x)=y.$ 
\qed

   \section{How to decompose a symmetric deleted join as a diagram}
\label{sec4}

 Let $P$ be a poset with $\hat 0$. We call two elements $x,y\in P$ {\bf disjoint}
if $P_{\leq x}\cap P_{\leq y}=\{\hat 0\}$, in other words from $z\leq x$, $z\leq y$ 
it follows that $z=\hat 0$. Denote $\hat P=P\setminus\{\hat 0\}$.

\begin{df}
 Let $A$ be a strictly increasing $t$-tuple of natural numbers, i.e. $A=(a_1,\dots,a_t)$,
where $a_1<\dots <a_t$, $a_i\in {\Bbb N}$. The {\bf deleted join} of $P$ is a poset 
$P^{[A]}$ consisting of all $t$-tuples $(x_1,\dots,x_t)$ with $x_i\in P$, such that 
$x_i$ and $x_j$ are disjoint for $1\leq i<j\leq t$. The order relations are given by 
$(y_1,\dots,y_t)<(x_1,\dots,x_t)$ iff $y_i<x_i$ for all $i\in[t]$. 
The terms of the deleted join are in natural bijection with elements of $A$.
When $A=\{1,\dots,n\}$ we simply write $P^{[n]}$ instead of $P^{[A]}$.
\end{df}

  When the poset $P$ is the face lattice of some simplicial complex, our definition
coincides with the one for the deleted join of a simplicial complex given by K.~Sarkaria
in \cite{Sa}. In the example we compute in the subsequent sections, $P$
will be a (possibly infinite) Boolean algebra, i.e. a face lattice of a simplex.

\begin{df}
  The {\bf deleted product} of $P$ is a poset $\hat P^{[A]}$, which is obtained from 
$P^{[A]}$ by imposing the additional condition that in every $t$-tuple $(x_1,\dots,x_t)$ 
we have $x_i\neq\hat 0$, for $i\in[t]$. Again $\hat P^{[\{1,\dots,n\}]}=\hat P^{[n]}$. 
\end{df}

 Let $\hat\bn=\bn\setminus\{\hat0\}$. For a poset $P$ with $\hat0$ the deleted join 
$P^{[n]}$ can be decomposed as a diagram $\cd_n$ of posets over $\tbn$ with posets 
assigned to different elements of $\bn$ in the following way: to $A\subseteq[n]$ we 
assign $\cd_n(A)=\hat P^{[A]}$, to every element of $A$ there is exactly one 
corresponding term in the product. If $B\subset A$, the map $p_{AB}:\cd_n(A)\lra\cd_n(B)$ 
is a projection of $\hat P^{[A]}$ onto $\hat P^{[B]}$ which is given by forgetting the 
terms which corresponds to elements in $A\setminus B$. 

  Let us show that the diagram $\cd_n$ is really a decomposition of $P^{[n]}$.

\begin{prop}
  For any poset $P$ and positive integer $n$, $\dl\cd_n = P^{[n]}$. 
\end{prop}

\pr Recall that the elements of $\dl\cd_n$ have the form $(A,e)$, 
where $A=\{a_1,\dots,a_t\}\subseteq[n]$ 
and $e=(x_{a_1},\dots,x_{a_t})\in\cd_n(A)=\hat P^{[A]}$. We define 
a map on the elements of the posets $\phi:\dl\cd_n\ra P^{[n]}$ by $\phi(A,(x_{a_1},\dots,
x_{a_t}))=(y_1,\dots,y_n)$, where $y_{a_i}=x_{a_i}$ for $i\in[t]$ and $y_j=\hat 0$ 
otherwise. It is obvious that $\phi$ is a bijection on the sets of elements of $\dl\cd_n$ 
and $P^{[n]}$.
Furthermore, $\phi$ is clearly order-preserving and hence a poset isomorphism.
\qed
  
  Let $\dbdn$, be the barycentric subdivision of $\cd_n$. By Definition \ref{bddf} 
the underlying poset of $\dbdn$ is $\bd(\tbn)$ and for every 
$x=(S_1\subset S_2\subset\dots\subset S_t)\in\bd(\tbn)$ 
we have $\dbdn(x)=\hat P^{[S_t]}$. 
If $y=(S_{i_1}\subset S_{i_2}\subset\dots\subset S_{i_m})<x$ 
for some $1\leq i_1<i_2<\dots<i_m\leq t$, according to the same definition, the map 
$f_{xy}:\dbdn(x)\rightarrow\dbdn(y)$ is defined by
$$f_{xy}=\begin{cases}
\text{ id, } & \text{ if } S_{i_m}=S_t,\\
\text{ projection map } \hat P^{[S_t]} \hookrightarrow \hat P^{[S_{i_m}]} ,& \text{ if } 
S_{i_1}\supset S_1.
\end{cases}$$   
  
  Now let us define an action of the symmetric group $\cs_n$ on the diagram $\dbdn$. 
Fix $\pi\in\cs_n$. As definitions \ref{dddf} and \ref{actp} suggest, we shall describe a
corresponding pair $(m,\tau)\in\cm_{\dd K}(\dbdn,\dbdn)$. $m\in\cf(\bd(\tbn),\bd(\tbn))$ 
is defined by
$$m:(S_1\subset S_2\subset\dots\subset S_t) \lra
(\pi(S_1)\subset\pi(S_2)\subset\dots\subset\pi(S_t)).$$ 
Further, we have to define $\tau_x:\dbdn(x)\ra\dbdn(m(x))$. Unwinding definitions
gives us $\dbdn(x)=\hat P^{[S_t]}$ and $\dbdn(m(x))=\hat P^{[\pi(S_t)]}$. The map 
$\pi:S_t\ra\pi(S_t)$ induces poset isomorphism $\tau_x:\hat P^{[S_t]}\ra
\hat P^{[\pi(S_t)]}$, which is exactly the map that we are looking for.

  Now, as we have an action of $\cs_n$ on $\dbdn$ we can consider a new diagram
$\tilde\cd_n=\dbdn/\cs_n$. The underlying poset of $\tilde\cd_n$ is (by coincidence) 
again $\tbn$, namely it is Bd$(\tbn)/\cs_n$. For $x=\{a_1,\dots,a_t\}\in\tbn$ 
with $a_1<\dots<a_t$ square (I) of \ref{cthm} asserts in our case that 
$$\tilde\cd_n(x)=\hat P^{[a_t]}/\cs_{a_1}\times\cs_{a_2-a_1}\times\dots\times\cs_{a_t-a_{t-1}},$$
where $\cs_{a_1}$ acts on the first $a_1$ terms of $\hat P^{[a_t]}$, $\cs_{a_2-a_1}$
acts on the next $a_2-a_1$ terms and so on.

       \section{An application: symmetric deleted join of a simplex}

  Let us now specify the choice of the partially ordered set $P$. Take $P=\cb_\infty$. 
Following the terminology of Section \ref{sec6}, $P^{[n]}=\cb_\infty^{[n]}=
\cb_{\infty,n}^o$. According to Theorem \ref{pmain}(a), $\da(\cb_{\infty,n}^o)$ 
is an infinite contractible simplicial complex, hence the spaces of the diagram 
$\da(\ti\cd_n)$ are given by  
\begin{multline}\label{cdnx}
\Delta(\tilde\cd_n(X))=\Delta(P^{[a_m]}/\cs_{a_1}\times\cs_{a_2-a_1}\times\dots
\times\cs_{a_m-a_{m-1}})= \\
=\Delta(P^{[a_m]})/\cs_{a_1}\times\cs_{a_2-a_1}\times\dots\times \cs_{a_m-a_{m-1}}
\simeq K(\cs_{a_1}\times\dots\times\cs_{a_m-a_{m-1}},1)\simeq \\
\simeq K(\cs_{a_1},1)\times\dots\times K(\cs_{a_m-a_{m-1}},1),
\end{multline}
where $X=\{a_1,\dots,a_m\}$.


  Let us understand the structure of the morphisms of the diagram $\da(\ti\cd_n)$.
For $x\geq y\in\hat\cb_n$, denote $X=\dbd_n(x)$, $Y=\dbd_n(y)$. Both $X$ and $Y$ 
are contractible and the diagram morphism $f:X\ra Y$ is $G$-equivariant and therefore 
induces a map $f/G:X/G\ra Y/G$. We have two fibre bundles $H\ra X\ra X/G$ and 
$K\ra Y\ra Y/G$, where $H$ and $K$ are subgroups of $G$, such that either $H\subseteq K$
or $K\subseteq H$. Since the map $f$ is $G$-equivariant, it induces a map between the two 
long exact sequences of homotopy groups associated to the fibre bundles. It follows from 
the naturality of these long exact sequences that the map $f/G$ is induced by the 
homomorphism $\iota:H\ra K$, where $\iota$ is inclusion if $H\subseteq K$ and restriction
if $K\subseteq H$. 

  Thus the diagram $\tilde\cd_n$ consists of classifying spaces for various
direct products of symmetric groups glued together by the maps that are induced by the
homomorphisms between these direct products.

  In the remainder of this section we compute $H_*(\da(\dl\tilde\cd_n),\zp)$. 
All homology groups will have coefficients in $\zp$, unless the contrary 
is explicitly stated. According to the Proposition \ref{2.4}, 
$H_*(\da(\dl\tilde\cd_n))=H_*(\holim(\Delta(\tilde\cd_n)))$. 
To calculate the homology groups of $\holim
(\Delta(\tilde\cd_n))$, we set up a spectral sequence $(E_r^{*,*})_{r=1}^\infty$,
using the following filtration:
$$F_0=\emptyset,\,\,\,F_m=\holim\tilde\cd_n|_{\tbn(m)},$$ 
where $\tbn(m)$ is the $m$-skeleton of $\tbn$ a poset which consists of all nonempty subsets of 
$[n]$ of cardinality at most $m$. Clearly $\emptyset=F_0\subset F_1\subset\dots
\subset F_n=\holim\tcdn$.

  To understand the differentials in this spectral sequence, let us describe a cell 
structure on $\holim\tcdn$. Let $X\subset[n]$. Since every $\tcdn(X)$ is a classifying 
space of a direct product of some symmetric groups, it can be represented by a simplicial 
complex (the particular structure of that complex is not important in our context). The 
building blocks of $\holim\tcdn$ are spaces of the form $\da((\tbn)_{\leq X})\times\tcdn(X)$. 
The space $T_X=\Delta((\tbn)_{\leq X})$ can be viewed as a simplex, i.e. a simplicial
complex having exactly one largest cell $T$. The cells of $T_X\times\tcdn(X)$ are 
direct products of the cells from the two simplicial complexes. Let $c$ be a cell 
of $\tcdn(X)$ and $t$ be a cell of $T_X$, then 
\begin{equation}\label{gd}
\bo(t\times c)=\bo t\times c+t\times\bo c.
\end{equation}

  Let us describe the entries in $E_1^{*,*}$. First of all we have
\begin{equation}\label{eq1}
H_*(F_m,F_{m-1})=\bigoplus_{X\subseteq[n],|X|=m}H_*(T_X\times\tcdn(X),
\delta T_X\times\tcdn(X)),
\end{equation}
where $\delta T_X$ is the boundary of the simplex $T_X$, i.e. $\delta T_X=T_X\sm T$. 
From \eqref{gd} and \eqref{eq1} we see that relative to $F_{m-1}$ we have
$\bo(t\times c)=0$ unless $t$ is the unique largest cell $T$ of $T_X$, in which case
$\bo(T\times c)=T\times\bo c$. This means that the boundary operator in 
$C_*(T_X\times\tcdn(X),\delta T_X\times\tcdn(X))$ is identical to the one in 
$C_*(\tcdn(X))$ up to the shift by the dimension of the simplex $T_X$, that is
$Z_*(T_X\times \tcdn(X),\delta T_X)=\{T\times c\,|\,c\in Z_*(\tcdn)\}$ and the same 
for $B_*$. Hence
\begin{equation}\label{e1**}
E_1^{m,d-m}=H_d(F_m,F_{m-1})=\bigoplus_{X\subseteq[n],|X|=m}H_{d+m-1}(\tcdn(X)).
\end{equation}

  Furthermore, in $\holim\tcdn$, whenever $t$ is not the largest cell of $T_X$, say
$t$ corresponds to some $Y\subset X$, the cell $t\times c$ gets identified with 
$t\times f_{X,Y}(c)$. In order to describe $d_1$ - the first differential of our 
spectral sequence we introduce an auxiliary diagram $\crxc$ for any 
$X\subseteq [n]$ and $c\in Z_*(\tcdn(X))$. The underlying poset of $\crxc$
is $(\tbn)_{<X}$. For $Y\subset X$, we have $\crxc(Y)=f_{X,Y}(c)$. For 
$Y_2\subset Y_1\subset X$, the diagram morphism between $\crxc(Y_1)$ and
$\crxc(Y_2)$ is induced by $f_{Y_1,Y_2}$. We define $f(c)=\holim\crxc$
and $\tilde f(c)=\holim(\crxc|_{P_x})$, where $P_x=(\tbn)_{<X}\sm\{$coatoms$\}$.

Let $X=\{a_1,\dots a_m\}$ and $c\in Z_*(\tcdn(X))$, $d_1$ is given by
\begin{equation}\label{d1}
d_1(T\times c)=(\delta T_X\times f(c),\tilde\delta T_X\times\tilde f(c))=
\sum_{i=m}^t (-1)^i\,T_i\times f_i(c),
\end{equation} 
where $f_i=f_{X,X\sm\{a_i\}}$, $T_i$ is a face corresponding to $X\sm\{a_i\}$, 
and $\tilde\delta T_X=\delta T_X\sm\{T_1,\dots,T_m\}$.

  We can immediately show that the spectral sequence collapses after the second
step. Namely, take such $c\in Z_*(\tcdn(X))$ that $d_1(c)=0$, this means that
$f_i(c)=0$ for $i\in [m]$. Let $Y\subset X$ and pick $i\in X\setminus Y$, then
$f_{X,Y}=f_{X\setminus\{i\},Y}\circ f_i=0$. This means that $\crxc(Y)=0$ for 
any $Y\subset X$ and hence $f(c)=\holim R_{X,c}=0$, therefore $\bo(c)=T\times\bo c+
\bo T\times f(c)=0$, which means that the sequence collapses.

In order to understand the maps $f_i$ on the homology level better we should first 
understand the homology of the symmetric group.  For that we give an 
excerpt from the work of M.~Nakaoka, \cite{N60,N61}, describing the homology groups 
$H_*(K(\cs_n,1),\zp)=H_*(S_n)$.

\begin{thm}\label{decomp} (\cite[Corollary 5.9]{N60})

 For any group $G$, $n$ and $q$ positive integers, we have
$$H_q(\cs_n,G)= H_q(\cs_{n-1},G)\oplus H_q(\cs_n,\cs_{n-1},G).$$
\end{thm}

The following is a comprised version of sections 6.1 and 6.2 of \cite{N60}. 

  Given a prime $p$, we denote by $Q(p)$ the set of all sequences
$J=(j_1,\dots,j_l)$, $l>0$, of positive integers satisfying the following
conditions:
\be
\item $j_k=0$ or $-1$ mod $2(p-1)$ for $1\leq k\leq l$,
\item $j_k\leq p j_{k+1}$ for $1\leq k<l$,
\item $j_1>(p-1)(j_2+\dots+j_l)$.
\ee  

  For each element $J\in Q(p)$, define the dimension and the rank by
$$D(J)=j_1+\dots+j_l,\,\,\,R(J)=p^l.$$

  Let $Q(P)_r^d$ denote the set of all elements $a\in Q(P)$, such that
$D(a)=d$ and $R(a)=r$, obviously $Q(P)=\cup_{d,r}Q(P)_r^d$.
Let $U(Q(P))$ be the $\zp$-algebra generated by all elements $a\in Q(P)$
subject to the relations 
\begin{equation} \label{antic}
ab=(-1)^{D(a)D(b)}ba.
\end{equation}
In particular, if $D(a)$ is odd and $p\neq 2$, then $a^2=0$.

We linearly extend functions $D$ and $R$ to all the monomials in $U(Q(P))$: 
\begin{equation} \label{dr}
D(ab)=D(a)+D(b),\,\,\,R(ab)=R(a)+R(b),
\end{equation}
for monomials $a,b\in U(Q(P))$. Again, for a monomial $\phi\in U(Q(P))$ we call $D(\phi)$ and $R(\phi)$ 
the dimension and the rank of $\phi$. We denote by $U_r^d(Q(P))$ the submodule 
generated by all monomials of dimension $d$ and rank $r$.

\begin{thm} \label{hom}(\cite[Theorem 6.3]{N60}).
  $H_d(\cs_m,\cs_{m-1};\zp)=U_m^d(Q(p))$, hence, because of the Theorem 
\ref{decomp}, we have $H_d(\cs_m,\zp)=\bigoplus_{r\leq m}U_r^d(Q(p))$.
\end{thm}
  
Theorem \ref{hom} can be rephrased as follows: the Betti number $\beta_d(\cs_m,\cs_{m-1})$ 
is equal to the number of choices $(J_1,\dots,J_t)$, $J_i\in Q(P)$ such that
\be
\item[(1)] $\sum_{i=1}^t R(J_i)=m$;
\item[(2)] if $J_i=J_j$ for some $i\neq j$ then either $D(J_i)$ is even or $p=2$ (we call this the 
{\it nonvanishing condition}).
\ee 

  Let $k$ and $m$ be positive integers. Naturally defined group homomorphisms 
$\iota:\cs_k\lra\cs_{k+m}$ and $\mu:\cs_k\times\cs_m\lra\cs_{k+m}$ (for $\pi_1\in\cs_k$, 
$\pi_2\in\cs_m$, $\mu(\pi_1\times\pi_2)$ acts as $\pi_1$ on the first $k$ elements of 
$[k+m]$ and as $\pi_2$ on the last $m$ elements of $[k+m]$; $\iota(\pi_1)=\mu(\pi_1\times
\text{id})$) induce homology maps $\iota_*:H_*(\cs_k)\lra H_*(\cs_{k+m})$ and $\mu_*:
H_*(\cs_k\times\cs_m)\lra H_*(\cs_{k+m})$. Nakaoka has combinatorially described these
maps: $\iota_*$ is a monomorphism (see \cite[Theorem 5.8]{N60} and description 
following it) and $\mu_*$ coincides with the multiplication operation in the algebra
$U(Q(P))$ (see \cite[Lemma 1.2,(4)]{N61}).

   Now, using these results, we are ready to resume the computation.

  Since the spectral sequence collapses, $H_*(\dl\ti\cd_n)=H_*(E_1^{*,*},d_1)$,
where $(E_1^{*,*},d_1)$ is the chain complex described by \eqref{e1**} and \eqref{d1}.
We denote this complex  by $\dan$ and shall now describe its combinatorial structure.

  The vector space $C_*(\dan)$ has a basis $\cb$ which can be indexed by triples
$(J,\pi,f)$, where $J=\{J_1,\dots,J_t\}$, $J_i\in Q(P)$; $\pi=(\pi_1,\dots,\pi_m)$; 
$f:[t]\lra[m]$ satisfying the following conditions:

{\bf Conditions A.}\nl
(1) $\,$ $\sum_{i=1}^t R(J_i)\leq\sum_{r=1}^m \pi_i\leq n$;\nl
(2) $\,$ $\sum_{f(i)=j}R(J_i)\leq\pi_j$, for all $j\in[m]$;\nl
(3) $\,$ the set $\{J_i\}_{f(i)=j}$ satisfies the nonvanishing condition for all $j\in[m]$.

  Combinatorially, conditions (2) and (3) mean that we distribute $J_i$'s into
blocks with sizes given by $\pi$ such that the total rank of $J_i$'s which are put
into the same block does not exceed its size and no two equal $J_i$'s of odd rank
are put into the same block (putting in the same block simply means taking product).
 
  Let $d=\sum_{i=1}^t D(J_i)$ then $(J,\pi,f)\in C_{d+m-1}(\dan)$. It is easy to read 
off the boundary map of $\dan$ from the description above:
\begin{equation}\label{ddf}
  \bo(J,\pi,f)=\sum_{i=1}^{m-1}(-1)^i(J,\pi^i,f^i)+(-1)^m\epsilon,
\end{equation} 
 where $\pi^i=(\pi_1,\dots,\pi_i+\pi_{i+1},\dots,\pi_m)$ and
$$f^i(j)=
\begin{cases}
  f(j),   & \text{ if } f(j)\leq i \\
  f(j)-1, & \text{ if } f(j)\geq i+1.
\end{cases}$$
  The $\epsilon$ above is given by the rule:
\begin{itemize}
\item if im$\,f\subseteq[m-1]$ then $\epsilon=(\pi',f')$, where 
$\pi'=(\pi_1,\dots,\pi_{m-1})$ and $f'=f|_{[m-1]}$;
\item if $m\in\,$im$\,f$, then $\epsilon=0$.
\end{itemize}

  It is easy to see that sets $\pi=\{\pi_1,\dots,\pi_m\}$, satisfying condition {\bf A}(1) 
above, are in bijection with nonempty subsets of $[n]$ (simply define $\{a_1,\dots,a_m\}
\subseteq[n]$ by $a_i=\sum_{j=1}^i\pi_j$). Thus the basis $\cb$ can be split into groups
indexed by nonempty subsets of $[n]$. This is exactly the picture we see on the first
step of our spectral sequence above, compare \eqref{cdnx} and \eqref{e1**}.

  To compute the second, terminal step of the spectral sequence we shall regroup
our basis elements. The key observation is that the boundary operator $\bo$ does
not change $J$, see \eqref{ddf}. Hence, if we fix the set $J=(J_1,\dots,J_t)$,
the chain complex $\dan(J)=\text{span}((J',\pi,f)\,|\,(J',\pi,f)\in\cb,\,J'=J)$ is a 
subcomplex of $\dan$. Its basis elements can be simply indexed by pairs $(\pi,f)$ (of 
course satisfying conditions {\bf A}). In what follows let $l_i=|J_i|$ and $r_i=p^{l_i}$.

  Clearly, the homology of $\da_n$ splits as a direct sum 
\begin{equation}\label{dsum}
   H_k(\da_n)=\bigoplus H_{k-d}(\da_n(J_1,\dots,J_t))
\end{equation}
where $r=\sum_{i=1}^t r_i$ and $d=\sum_{i=1}^q D(J_i)$ and the sum is taken over 
all sets $\{J_1,\dots,J_t\}$ satisfying $r\leq n$.

  It turns out that the parameters $l_i$ (or equivalently $r_i$) described above are 
excessive for our computation. Let us describe yet another chain complex which will be 
isomorphic to $\da_n(J_1,\dots,J_t)$. The chain complex that we shall 
now describe will be denoted $\da_n(\co,\ce,t)$, where $\co=\{\co_1,\dots,\co_u\}$,
$\ce=\{\ce_1,\dots,\ce_v\}$, such that $\co_i,\ce_j\subseteq[t]$, sets $\{ \co_i \}$ and
$\{ \ce_j \}$ are disjoint and $\sum_{i=1}^u|\co_i|+\sum_{j=1}^v|\ce_j|\leq t$. 

{\bf Description of $\da_n(\co,\ce,t)$.} The elements of the basis of $C_*(\da_n
(\co,\ce,t))$ are indexed by pairs $(\pi,f)$ with $\pi=(\pi_1,\dots,\pi_m)$ and $f:[t]\lra[m]$,
satisfying following conditions:

{\bf Conditions B.}

(1) $t\leq\sum_{i=1}^m\pi_i\leq n$;

(2) $|f^{-1}(j)|\leq\pi_j$ for all $j\in[m]$;

(3) Assume that $x<y$. If $x,y\in\co_i$ for some $i\in[u]$ then $f(x)<f(y)$.
If $x,y\in\ce_i$ for some $i\in[v]$ then $f(x)\leq f(y)$.    

 Clearly, the difference from the conditions {\bf A} is that the rank function $R$ has 
been replaced by constant 1 and the set $\co$, resp. $\ce$, keeps track of the sets of 
mutually equal elements of $J$ which have odd, resp. even, dimension. As before $(\pi,f)
\in C_{d+m-1}(\da(\co,\ce,t))$, where $d=\sum_{i=1}^t D(J_i)$. 

The differential in this new chain complex is defined in the same way as above, 
see \eqref{ddf}.

\begin{lm} \label{redl}
  The chain complexes $\da_n(J_1,\dots,J_t)$ and $\da_{n'}(\co,\ce,t)$ are isomorphic.
Here $n'=n+t-(r_1+\dots+r_t)$, $\co=\{\co_1,\dots,\co_u\}$ and $\ce=\{\ce_1,\dots,\ce_v\}$.
If $p\neq 2$, $\co_i$'s, resp. $\ce_i$'s, are all sets of elements of $J$ which are mutually 
equal and have odd, resp. even, dimension. If $p=2$, then $\co_i$'s are all sets of elements 
of $J$ which are mutually equal and $\ce=\emptyset$.
\end{lm} 

\pr An isomorphism $\phi$ is defined by $\phi(\pi,f)=(\pi',f)$, where 
$\pi'=(\pi_1-\alpha_1,\dots,\pi_t-\alpha_t)$, here $\alpha_j=\sum_{f(i)=j}(r_i-1)$.
\qed

\begin{prop} \label{dnt}
  $\dan(\co,\ce,t)$ is acyclic unless $t=n$ and $v=0$, i.e. $\ce=\emptyset$, in which
case $\dan(\co,\ce,t)$ has homology of an $(n-1)$-sphere, that is
  \begin{equation}\label{homn}
    \th_i(\da_n(\co,\emptyset,n))=
    \begin{cases}
        \Bbb Z, & \text{ if } i=n-1,\\
          0,    & \text{ otherwise. }
    \end{cases}
  \end{equation}
\end{prop}

\pr Let $\ca_{n+1}$ be the braid arrangement in ${\Bbb R}^{n+1}$ and let $\Lambda_{n+1}$
be the collection of cells of the stratification of the sphere $S^{n+1}$ by $\ca_{n+1}$.
We shall now describe how to associate a cell from $\Lambda_{n+1}$ to an element of the
basis of $\dan(\co,\ce,t)$.

  Let $(\pi,f)$ be such an element with $\pi=(\pi_1,\dots,\pi_m)$ and $f:[t]\ra[m]$.
First we add an artificial block $\pi_{m+1}$ to $\pi$ such that $|\pi_{m+1}|=n+1-t$.
The function $f$ shows how $J_i$'s are distributed in blocks $\pi_1,\dots,\pi_m$. Let us 
augment the set $\{J_1,\dots,J_t\}$ by $n-t$ elements called $*_{t+1},\dots,*_{n+1}$,
which we distribute into the blocks $\pi_1,\dots,\pi_{m+1}$ filling the places left over
from $J_i$'s. More formally, the function $\ti f:[n+1]\ra[m+1]$ is uniquely defined by

(1) $\ti f|_{[t]}=f$;

(2) $\ti f|_{\{t+1,\dots,n+1\}}$ is increasing and $|\ti f^{-1}(j)|=\pi_j$ for all 
$j\in[m+1]$. 

   The cell $c_{\pi,f}\in\Lambda_{n+1}$ associated to $(\pi,f)$ is given by the
equations on $x_1,\dots,x_{n+1}$ which are induced by the function $\ti f$ in the
obvious manner: if $\ti f(i)=\ti f(j)$ then $x_i=x_j$; if $\ti f(i)<\ti f(j)$ then
$x_i<x_j$.

   Let $\cc=\cup_{(\pi,f)}c_{\pi,f}$, where the union is taken over whole basis
of $\dan(\co,\ce,t)$. It is not difficult to see that the boundary operator of the
chain complex $\dan(\co,\ce,t)$ coincides with the natural boundary operator
in $\Lambda_{n+1}$, hence $\dan(\co,\ce,t)\simeq(\ca_{n+1},\ca_{n+1}\sm\cc)$
(considering $\ca_{n+1}$ relative to $\ca_{n+1}\sm\cc$ means exactly that we limit
our attention to the cells in $\cc$). The whole set $\cc$ can be described by the
following inequalities:

{\bf Inequalities C.}

(1) $x_i<x_{n+1}$ if $i\in[t]$ (since we have to make sure that $n+1$ is in the block
$\pi_{m+1}$ and none of the $J_i$'s is);

(2) $x_i\leq x_j$ if $t+1\leq i<j\leq n+1$;

(3) $x_i<x_j$ if $i,j\in\co_k$ for some $k\in[u]$;

(4) $x_i\leq x_j$ if $i,j\in\ce_k$ for some $k\in[v]$.

  Clearly, if all of the inequalities were non-strict, then $\cc$ would be a closed convex
set, homeomorphic to an $(n-1)$-ball $T$. However, even in our case it is clear that
$(\ca_{n+1},\ca_{n+1}\sm\cc)\simeq(T,Z)$, where $Z$ is a piece of the boundary of that 
ball. If $n=t$ and $\ce=\emptyset$ then we only have strict inequalities, hence $Z$
is the full boundary of $T$ and $(T,Z)\simeq (S^{n-1},*)$ where $*$ is a single point, 
which proves \eqref{homn}.

  Let us now show that otherwise $Z$ is contractible (that would mean that $(T,Z)$
is acyclic and hence so is $\dan(\co,\ce,t)$). We call an inequality from $\bf C$ 
{\it essential} if it does not follow from other inequalities from $\bf C$. $Z$ can be
seen as a union of convex sets $(Z_i)_{i\in I}$, each $Z_i$ is defined by changing
one of the essential inequalities in (1) or (3) to an equality and making other
inequalities non-strict. Let us show that all of these convex sets have a common point.

   If $t\geq 1$, let $\alpha$ be the cell defined by $x_{t+1}\leq\dots\leq x_{n+1}=
x_1=\dots=x_t$, clearly $\alpha\in Z_i$, for all $i\in I$, since all of the strict
inequalities (which are turned to equalities in $Z_i$) contain only elements from
the set $\{1,\dots,t,n+1\}$.

  Assume now that $t=n$ but $\ce\neq\emptyset$, say $c,d\in\ce_1$, $c<d$. Let 
$\alpha$ be defined by $x_c\leq x_1=\dots=x_{c-1}=x_{c+1}=\dots=x_{n+1}$. The only 
strict inequality in $\bf C$ containing $x_c$ is $x_c<x_{n+1}$, which is a consequence
of $x_c\leq x_d$ and $x_d<x_{n+1}$, hence it is not an essential inequality. 
Therefore we may conclude that also in this case $\alpha\in Z_i$ for all $i\in I$.
  
   Since in both cases $\cap_{i\in I}Z_i$ and $Z_i$'s are convex, $Z=\cup_{i\in I}Z_i$
is contractible. This finishes the proof of the proposition.
\qed

 Let us now proceed to the Main Theorem. Recall Nakaoka's description of the relative homology
$H_*(\cs_n,\cs_{n-1};\zp)$. To obtain $H_*(\dl\tcdn,\zp)$ we only need to make one 
modification in his description:
the dimension function should be different: $\ti D(a)=D(a)+1$ for all $a\in Q(P)$.

Let the algebra with this new dimension function be denoted by $\tilde U(Q(P))$, then
\begin{mthm}\label{main}
\begin{equation}\label{emain}
  H_q(\da(\dl\tcdn),\zp)=\tilde U_n^{q+1}(Q(p)).
\end{equation}
\end{mthm}

\pr According to \eqref{dsum} we have
\begin{equation}\label{e1}
H_q(\da_n)=\bigoplus H_{q-d}(\da_n(J_1,\dots,J_t)),
\end{equation}
where $d=\sum_{i=1}^t D(J_i)$, $r=\sum_{i=1}^t R(J_i)$ and the sum is taken over all 
$t$-tuples $(J_1,\dots,J_t)$ satisfying $r\leq n$. Furthermore, Lemma \ref{redl}
asserts that 
\begin{equation}\label{e2}
H_{q-d}(\da_n(J_1,\dots,J_t))=H_{q-d}(\da_{n'}(\co,\ce,t)),
\end{equation}
where $n'=n+t-r$. 

Combining \eqref{e1}, \eqref{e2} and Proposition \ref{dnt} we can conclude that 
$\beta_q(\da_n)$ is equal to the number of all $t$-tuples $(J_1,\dots,J_t)$ such that 
$\ce=\emptyset$, $r=n$ and $d=q+1-t$, i.e. $\ti d=q$, where $\ti d=\sum_{i=1}^t\ti D(J_i)$, 
which proves this theorem.
\qed

      \section{Homology of the $k$-fold Boolean algebras}
\label{sec6}

 Throughout this section $n$ is either a positive integer or $\infty$.
\begin{df} \label{dfbnk}
  Let $k$ be a positive integer and $n$ either a positive integer or $\infty$.
We define posets $\bnk$, which we call {\bf $k$-fold Boolean algebras}, as follows. 
Elements of $\bnk$ are sets, each consisting of $k$ disjoint, possibly empty, finite 
subsets of $[n]$ (here $[\infty]=\Bbb N$), in other words, $\{S_1,\dots,S_k\}$ such 
that $S_i\subseteq[n]$, $S_i\cap S_j=\emptyset$ for $i\neq j$ and $|S_i|<\infty$. 
The partial ordering is given by the rule: $\{S_1',\dots,S_k'\}\leq\{S_1,\dots,S_k\}$ 
if and only if there exists a permutation $\pi$ of the set $[k]$ such that 
$S_i'\subseteq S_{\pi_i}$ for all $i\in[k]$.  
\end{df}
  
Using the terminology of Section \ref{sec4}, let us define $\bnko=\bn^{[k]}$ and
$\hbnko=\hat\bn^{[k]}$. Then the Definition \ref{dfbnk} can be reformulated as
$\bnk=\bnko/\cs_k$. Analogously, we define $\hbnk=\hbnko/\cs_k$. The latter can
be directly defined just like $\bnk$ with an extra condition that the sets $S_i$
should be non-empty.

  The following theorem is the main result of this section.
\begin{thm}\label{pmain} $\,$

(a) $\da(\bnko)$ is homotopy equivalent to a wedge of $(n-1)$-dimensional spheres
if $n<\infty$ and $\da(\cb_{\infty,k}^o)$ is contractible;
  
(b) homology groups of $\da(\bnk)$ with integral coefficients converge to those
of $\da(\cb_{\infty,k})$, more precisely
\begin{equation}
\ti H_i(\bnk,\Bbb Z)=
\begin{cases}
  \th_i(\cb_{\infty,k},\Bbb Z) & \text{ for } 0\leq i\leq n-2,\\
  {\Bbb Z}^{\mu(\bnk)} & \text{ for } i=n-1,\\
  0 & \text{ otherwise;}
\end{cases}
\end{equation}

(c) $H_q(\cb_{\infty,k},\zp)=\ti U_{k}^{q+1}(Q(p))$;

(d) same as (a) with $\bnko$ and $\cb_{\infty,k}^o$ replaced by $\hbnko$ and 
$\hat\cb_{\infty,k}^o$, and $(n-1)$-dimensional replaced by $(n-k)$-dimensional;

(e) same as (b) with $\bnk$ and $\cb_{\infty,k}$ replaced by $\hbnk$ and $\hat\cb_{\infty,k}$
, and $(n-1)$-dimensional replaced by $(n-k)$-dimensional.

\end{thm}

  We need some preparation before we can proceed with the proof of Theorem \ref{pmain}.
The following generalizes the definitions of $\bnk$ and $\bnko$.  

\begin{df} \label{dfbnkt}
  Let $t$ and $k$ be integers, $0\leq t\leq k$ and $n\in [\infty]$. An element of the 
poset $\bnkt$ is a pair consisting of a $t$-tuple $(S_1,\dots,S_t)$ (the ordered part) 
and a set $\{S_{t+1},\dots,S_k\}$ (the unordered part) such that $S_i\subseteq [n]$,
$S_i\cap S_j=\emptyset$ for $i\neq j$ and $S_i$ is finite. The partial ordering is given 
by the rule: $$(S_1',\dots,S_t')\{S_{t+1}',\dots,S_k'\}\leq(S_1,\dots,S_t)\{S_{t+1},\dots,
S_k\}$$ if and only if $S_i'\subseteq S_i$ for $i=1,\dots,t$ and there exists a permutation 
$\pi$ of the set $\{t+1,\dots,k\}$ such that $S_i'\subseteq S_{\pi_i}$ for $i=t+1,\dots,k$.   
\end{df}

  Clearly ${\cal B}_{n,k,0}=\bnk$ and ${\cal B}_{n,k,n}=\cb_{n,k,n-1}=\bnko$.

Observe that the lower intervals of $\bnkt$ are boolean algebras, hence every
$\bnkt$ is a face poset of a regular CW complex.

\begin{lm} \label{bnktsh}
  If $t\geq 1$ then the poset $\bnkt$ is EL-shellable, see \cite{BW83}, in particular
$\Delta(\bnkt)$ is homotopy equivalent to a wedge of spheres of dimension $n-1$.
\end{lm}

\pr Assume that $t\geq 1$ and $k\geq 2$. Our labeling will use integer labels.
Observe that every edge $(x,y)$ in the Hasse diagram of 
$\bnkt$ is either of the form $(x,\hat 1)$ or a merging of one of the sets of $x$ with
some element $a\in [n]$. Since $t>0$ we have a certain set in $x$ which is first in
the ordering of the corresponding $t$-tuple. In the rest of the proof we often refer to 
the sets of $x$ or $y$ as blocks. 

We divide all edges in the Hasse diagram 
into three groups and impose the labels in every group separately:
\be
\item[(1)] edges $(x,y)$, where $y\neq\hat 1$, and the element $a$ described above
is merged to the first block of $x$, we label such edge with $-a$; 
\item[(2)] edges from some coatom $x$ to $\hat 1$, we label these edges with $0$;
\item[(3)] edges $(x,y)$, where $y\neq\hat 1$, and the element $a$ described above
is merged to a block of $x$ which is not the first one, we label such edge with $a$. 
\ee
 Let us now show that these labels give us an EL-labeling. Consider an interval $(x,y)$.  
We divide the further proof into two cases.

{\bf Case 1.} $y\neq\hat 1$. The increasing chain is given by merging first all elements
which have to be merged to the first block in decreasing order and then all other
elements in increasing order. This chain is obviously lexicographically least.
It is also uniquely defined, since one has to merge the elements to the first block
first in the unique way, and then there is a unique way to merge the rest 
of the elements to the other blocks (creating a new nonempty block is considered as a 
merging to an empty block).

{\bf Case 2.} $y=\hat 1$. The last label in any maximal chain in $(x,y)$ has the label $0$,
hence every increasing chain should have only nonpositive labels. This means that we have 
to merge all the elements which are not already in one of the blocks of $x$ to the first
block, and we have to do it in decreasing order. Hence again the increasing chain is
unique and lexicographically least. 
  
 We have shown that for $t>0$ the posets $\bnkt$ are EL-shellable, hence $\bnkt$
are shellable and the simplicial complexes $\Delta(\bnkt)$ are homotopy equivalent 
to a wedge of spheres.
\qed

  The following technical lemma is needed for the transition from the sequence of
finite posets to their infinite limit.

\begin{lm} \label{apprlm}
  Assume that we have a family of finite simplicial complexes $\{K_i\}_{i=1}^\infty$
and an infinite simplicial complex $K_\infty$, such that $K_i$ is a subcomplex of $K_j$,
whenever $i<j$ or $j=\infty$, and $K_\infty=\cup_{i=1}^\infty K_i$. 

(1) If $K_i$ is $t_i$-connected and $\lim_{i\to\infty}t_i=\infty$ then $K_\infty$ 
is contractible.
 
(2) Assume that $\dim K_i=i-1$ and that for $t\geq 2$, $K_t$ is homotopy equivalent
to a regular CW complex obtained from $K_{t-1}$ by pasting on a number of $(t-1)$-cells.

   Then $\th_i(K_\infty)=\th_i(K_t)$, for all $i=0,1,2,\dots,t-2$.
\end{lm}

\pr (1) The image $A$ of a mapping of any compact space (in particular an $n$-sphere)
into $K_\infty$ lies in a union of finite number of cells. Hence for sufficiently large 
$N$, $A\subseteq K_N$ and $t_N\geq n$, which means that $A$ can be shrunk
to a point in $K_N$, hence also in $K_\infty$. Thus all homology groups of $K_\infty$
are trivial and therefore $K_\infty$ is contractible. 

(2) Obvious.
\qed

{\bf Proof of Theorem \ref{pmain}.} (a) Setting $t=n$ in Lemma \ref{bnktsh} we get 
that $\da(\bnko)$ is homotopy equivalent to a wedge of $(n-1)$-spheres, in particular
it is $(n-2)$-connected. Furthermore, Lemma \ref{apprlm} with $K_i=\da(\cb_{i,k}^o)$
asserts that $\da(\cb_{\infty,k}^o)$ is contractible.

(b) Recall that $\bnk$ is a face lattice of some regular CW complex.  Denote this
complex by $\ti\da_{n,k}$. Obviously $\ti\da_{n,k}\simeq\da(\bnk)$. Let us study how
the complex $\ti\da_{n+1,k}$ differs from $\ti\da_{n,k}$.

  It is clear that $\ti\da_{n+1,k}$ has a vertex $x$, represented by a set $(\{n+1\},
\emptyset,\dots,\emptyset)$, such that a cell of $\ti\da_{n+1,k}$ is not a cell of 
$\ti\da_{n,k}$ if and only if it contains $x$.  In other words $\ti\da_{n+1,k}=\ti
\da_{n,k}\cup\st_{\ti\da_{n,k}}(x)$. It is easy to see that $\st_{\ti\da_{n,k}}(x)\cap
\ti\da_{n,k}=\lk_{\ti\da_{n,k}}(x)\simeq\da(\cb_{n-1,k,1})$. The latter is, by Lemma 
\ref{bnktsh}, homotopy equivalent to a wedge of $(n-2)$-spheres. This means that $\ti
\da_{n+1,k}$ is homotopy equivalent to a space obtained from $\ti\da_{n,k}$ by pasting 
on a number of $(n-1)$-dimensional cells (one cell over each $(n-2)$-sphere in $\da
(\cb_{n-1,k,1})$). Thus the conditions of (2) of Lemma \ref{apprlm} are satisfied 
and hence $\th_i(\bnk,{\Bbb Z})=\th_i(\cb_{\infty,k},{\Bbb Z})$ for $0\leq i\leq n-2$.

  That $\th_{n-1}(\bnk,{\Bbb Z})={\Bbb Z}^{\mu(\bnk)}$ follows from the fact that 
$\da(\cb_{\infty,k})$ is $\Bbb Q$-acyclic.

(c) Clearly $\da(\cb_{\infty,k})\simeq\da(\dl\ti\cd_k)$, so Main Theorem \ref{main} applies.

(d) We shall assume that $n\geq k+2$, the other cases are obvious. $\bnko$ can be viewed 
as a face poset of a regular CW complex. Denote this complex by $C$. Let $P_i$ be the 
subposet of $\bnko$ consisting of elements $(S_1,\dots,S_k)$, such that $S_i=\emptyset$, 
where $i\in[k]$. Clearly, $P_i$ is a face poset of $C_i$, where $C_i$ is a regular CW 
subcomplex of $C$. Let $T=\cup_{i=1}^k C_i$. 

   There is a projection map $p$ from $C$ to a $k$-simplex $S$, defined by mapping a
$k$-tuple $(S_1,\dots,S_k)$ to its support subset of $[k]$, under which $T$ is mapped
to the boundary of $S$. $\da(\hbnk)$ is the preimage of the biggest cell of that simplex
and therefore $H_*((\da(\hbnk))\times(S,\delta S))=H_*(C,T)$. To show that $\da(\hbnk)$ is 
homotopy equivalent to a wedge of spheres, it is enough to show that $\pi_1(\da(\hbnk))=0$
and that $H_*(C,T)$ is zero except for the highest dimension.

  We know that $C_{i_1}\cap\dots\cap C_{i_t}$ is $(n-t-2)$-connected so by the 
Nerve Lemma, $T=\cup_{i=1}^k C_i$ is $(n-3)$-connected.  Thus $T$ has 
the homology of a wedge of spheres of dimension one less than that of $C$. Using
Mayer-Vietoris long exact sequence one sees that $H_*(C,T)$ (and hence $H_*(\da(\hbnk))$)
can be non-zero only in the highest dimension.

  Let us now show that $\pi_1(\da(\hbnk))=0$. Let $R$ be the rank selection of $\hbnk$ 
where we select the elements of rank at most 3. Assume that $\da(R)$ is simply connected.
We can add elements of $\hbnk\sm R$ in some linear extension order, to finally form the
full $\hbnk$. When we add an element $x\in\hbnk\sm R$, we are gluing a cone over
$\da((\hbnk)_{<x})$, which is connected, since it is a direct product of simplices.
Therefore, by repetitive use of Seifert - Van Kampen theorem, $\da(\hbnk)$ is simply
connected. 

  Let us now show that $\pi_1(R)=0$. $R$ is a face poset of a CW complex $W$. $W$ can
be described combinatorially as follows: vertices are ordered $k$-tuples of different
elements of $[n]$, edges are pairs of such ordered $k$-tuples differing in exactly
one coordinate, and the two-dimensional faces are all triangles and squares which
are formed by these edges. Assume that $\pi_1(R)\neq 0$ and take the shortest path 
which cannot be contracted. It must contain at least two changes of the same coordinate.
Moving the chosen path along squares we can permute the coordinates of any two consecutive 
changes.  Thus we can retract our path to a path of the same length where two consecutive
changes occur in the same coordinate. Now, we can use the proper triangle and replace
these two steps by just one step, combining these two changes into one change.
We thereby obtain a shorter path and hence a contradiction.  
 
(e) $\cs_k$ acts freely on both $\da(\hat\cb_{\infty,k}^o)$ and $\da(\hbnko)$. 
Furthermore, $\da(\hat\cb_{\infty,k}^o)$ is contractible and $\da(\hbnko)$ is 
$(n-k-1)$-connected, hence $\da(\hat\cb_{\infty,k})\simeq K(\cs_k,1)$ and 
$\th_i(\hbnko)=\th_i(K(\cs_k,1))$, for $i=0,\dots,n-k-1$.
\qed

\begin{thebibliography}{BW84l}

\bibitem[Bj84]{Bj84} A.~Bj\"orner, {\em Posets, regular CW complexes and 
Bruhat order}, Europ. J. Combinatorics {\bf 5}, (1984), pp. 7--16.

\bibitem[B\'e]{Be}  J.~B\'enabou, {\em Introduction to Bicategories}, Reports
of the Midwest Category Seminar, Lecture Notes in Mathematics {\bf 47},
(1967), pp. 1--77, Springer-Verlag.

\bibitem[BK]{BK} A.K.~Bousfield, D.M.~Kan, {\em Homotopy Limits, Completions
and Localizations, Part} II, Lecture Notes in Mathematics {\bf 304}, (1972),
Springer-Verlag. 

\bibitem[BW83]{BW83} A. Bj\"orner, M. Wachs, {\em  On lexicographically
shellable posets,} Trans. Amer. Math. Soc. {\bf 277}, (1983), pp. 323--341.

\bibitem[Da]{Da} V.I.~Danilov, {\em The geometry of toric varieties},
Russian Math. Surveys {\bf 33}(2), (1978), pp. 97--151.

\bibitem[Er]{Er} C.~Ehresmann, {\em Cat\'egories doubles et cat\'egories
structur\'ees}, C.R.~Acad. Sci., Paris, {\bf 256}, (1963), pp. 1198--1201.

\bibitem[FL]{FL} R.~Fritsch, D.M.~Latch, {\em Homotopy inverses for nerve},
Math. Z., {\bf 177}, (1981), pp. 147--179.

\bibitem[Ha]{Ha} P.~Hanlon, {\em A note on the homology of signed posets}, J. 
Alg. Comb. {\bf 5}, (1996), no. 3, pp. 245--250.

\bibitem[Ma]{Ma} S.~MacLane, {\em Categories for the Working mathematician},
Springer-Verlag, New York, 1971.

\bibitem[N60]{N60} M.~Nakaoka, {\em Decomposition theorem for homology
groups of symmetric groups}, Ann. Math., Vol. {\bf 71}, No. 1, (1960), 
pp. 17--42.

\bibitem[N61]{N61} M.~Nakaoka, {\em Homology of infinite symmetric group},
Ann. Math., Vol. {\bf 73}, No. 2, (1961), pp. 229--257.

\bibitem[Pa]{Pa} P.H.~Palmquist, {\em The double category of adjoint squares},
Reports of Midwest Category Seminar, Lecture Notes in Mathematics {\bf 195},
(1971), pp. 123--154, Springer-Verlag.

\bibitem[Q73]{Q73} D.~Quillen, {\em Higher algebraic K-theory} I, Lecture
Notes in Mathematics {\bf 341}, (1973), pp. 85--148, Springer-Verlag.

\bibitem[Q78]{Q78} D. Quillen, {\em Homotopy properties of the poset of
nontrivial $p$-subgroups of a group,} Adv. in Math. {\bf 28},
(1978), pp. 101--128.

\bibitem[Sa]{Sa} K.~Sarkaria, {\em A generalized Kneser conjecture}, 
J. Comb. Theory $B$, (to appear).

\bibitem[Se]{Se} G.~Segal, {\em Classifying spaces and spectral sequences},
Inst. Hautes. \'Etudes Sci. Publ. Math. No {\bf 34}, (1968), pp. 105--112. 

\bibitem[St]{St} R.P. Stanley, {\em Enumerative Combinatorics}, vol. I,
 Wadsworth, Belmont, CA, 1986.

\bibitem[WZ\v{Z}]{WZZ} V.~Welker, G.M.~Ziegler, R.~\v{Z}ivaljevi\'{c},
{\em Comparison lemmas and applications for diagrams of spaces},
preprint, 35 pages.

\bibitem[Z\v{Z}]{ZZ} G.M. Ziegler, R. \v{Z}ivaljevi\'{c}, {\em 
Homotopy types of subspace arrangements via diagrams of spaces,}
Math. Ann. {\bf 295}, (1993), pp. 527--548.

\end{thebibliography}
\end{document}